\documentstyle[epic,eepic,11pt,fullpage]{article}

\catcode`\@=11
\def\relaxnext@{\let\next\relax}
\font\tenmsx=msxm10
\font\sevenmsx=msxm7
\font\fivemsx=msxm5
\font\tenmsy=msym10
\font\sevenmsy=msym7
\font\fivemsy=msym5
\newfam\msxfam
\newfam\msyfam
\textfont\msxfam=\tenmsx  \scriptfont\msxfam=\sevenmsx
  \scriptscriptfont\msxfam=\fivemsx
\textfont\msyfam=\tenmsy  \scriptfont\msyfam=\sevenmsy
  \scriptscriptfont\msyfam=\fivemsy

\def\hexnumber@#1{\ifcase#1 0\or1\or2\or3\or4\or5\or6\or7\or8\or9\or
	A\or B\or C\or D\or E\or F\fi }

\font\teneuf=eufm10
\font\seveneuf=eufm7
\font\fiveeuf=eufm5
\newfam\euffam
\textfont\euffam=\teneuf
\scriptfont\euffam=\seveneuf
\scriptscriptfont\euffam=\fiveeuf
\def\frak{\relaxnext@\ifmmode\let\next\frak@\else
 \def\next{\errmessage{Use \string\frak\space only in math
mode}}\fi\next}
\def\goth{\relaxnext@\ifmmode\let\next\frak@\else
\def\next{\errmessage{Use \string\goth\space only in math
mode}}\fi\next}
\def\frak@#1{{\frak@@{#1}}}
\def\frak@@#1{\fam\euffam#1}

\edef\msx@{\hexnumber@\msxfam}
\edef\msy@{\hexnumber@\msyfam}

\mathchardef\boxdot="2\msx@00
\mathchardef\boxplus="2\msx@01
\mathchardef\boxtimes="2\msx@02
\mathchardef\square="0\msx@03
\mathchardef\blacksquare="0\msx@04
\mathchardef\centerdot="2\msx@05
\mathchardef\lozenge="0\msx@06
\mathchardef\blacklozenge="0\msx@07
\mathchardef\circlearrowright="3\msx@08
\mathchardef\circlearrowleft="3\msx@09
\mathchardef\rightleftharpoons="3\msx@0A
\mathchardef\leftrightharpoons="3\msx@0B
\mathchardef\boxminus="2\msx@0C
\mathchardef\Vdash="3\msx@0D
\mathchardef\Vvdash="3\msx@0E
\mathchardef\vDash="3\msx@0F
\mathchardef\twoheadrightarrow="3\msx@10
\mathchardef\twoheadleftarrow="3\msx@11
\mathchardef\leftleftarrows="3\msx@12
\mathchardef\rightrightarrows="3\msx@13
\mathchardef\upuparrows="3\msx@14
\mathchardef\downdownarrows="3\msx@15
\mathchardef\upharpoonright="3\msx@16

\mathchardef\downharpoonright="3\msx@17
\mathchardef\upharpoonleft="3\msx@18
\mathchardef\downharpoonleft="3\msx@19
\mathchardef\rightarrowtail="3\msx@1A
\mathchardef\leftarrowtail="3\msx@1B
\mathchardef\leftrightarrows="3\msx@1C
\mathchardef\rightleftarrows="3\msx@1D
\mathchardef\Lsh="3\msx@1E
\mathchardef\Rsh="3\msx@1F
\mathchardef\rightsquigarrow="3\msx@20
\mathchardef\leftrightsquigarrow="3\msx@21
\mathchardef\looparrowleft="3\msx@22
\mathchardef\looparrowright="3\msx@23
\mathchardef\circeq="3\msx@24
\mathchardef\succsim="3\msx@25
\mathchardef\gtrsim="3\msx@26
\mathchardef\gtrapprox="3\msx@27
\mathchardef\multimap="3\msx@28
\mathchardef\therefore="3\msx@29
\mathchardef\because="3\msx@2A
\mathchardef\doteqdot="3\msx@2B

\mathchardef\triangleq="3\msx@2C
\mathchardef\precsim="3\msx@2D
\mathchardef\lesssim="3\msx@2E
\mathchardef\lessapprox="3\msx@2F
\mathchardef\eqslantless="3\msx@30
\mathchardef\eqslantgtr="3\msx@31
\mathchardef\curlyeqprec="3\msx@32
\mathchardef\curlyeqsucc="3\msx@33
\mathchardef\preccurlyeq="3\msx@34
\mathchardef\leqq="3\msx@35
\mathchardef\leqslant="3\msx@36
\mathchardef\lessgtr="3\msx@37
\mathchardef\backprime="0\msx@38
\mathchardef\risingdotseq="3\msx@3A
\mathchardef\fallingdotseq="3\msx@3B
\mathchardef\succcurlyeq="3\msx@3C
\mathchardef\geqq="3\msx@3D
\mathchardef\geqslant="3\msx@3E
\mathchardef\gtrless="3\msx@3F
\mathchardef\sqsubset="3\msx@40
\mathchardef\sqsupset="3\msx@41
\mathchardef\vartriangleright="3\msx@42
\mathchardef\vartriangleleft="3\msx@43
\mathchardef\trianglerighteq="3\msx@44
\mathchardef\trianglelefteq="3\msx@45
\mathchardef\bigstar="0\msx@46
\mathchardef\between="3\msx@47
\mathchardef\blacktriangledown="0\msx@48
\mathchardef\blacktriangleright="3\msx@49
\mathchardef\blacktriangleleft="3\msx@4A
\mathchardef\vartriangle="0\msx@4D
\mathchardef\blacktriangle="0\msx@4E
\mathchardef\triangledown="0\msx@4F
\mathchardef\eqcirc="3\msx@50
\mathchardef\lesseqgtr="3\msx@51
\mathchardef\gtreqless="3\msx@52
\mathchardef\lesseqqgtr="3\msx@53
\mathchardef\gtreqqless="3\msx@54
\mathchardef\Rrightarrow="3\msx@56
\mathchardef\Lleftarrow="3\msx@57
\mathchardef\veebar="2\msx@59
\mathchardef\barwedge="2\msx@5A
\mathchardef\doublebarwedge="2\msx@5B
\mathchardef\angle="0\msx@5C
\mathchardef\measuredangle="0\msx@5D
\mathchardef\sphericalangle="0\msx@5E
\mathchardef\varpropto="3\msx@5F
\mathchardef\smallsmile="3\msx@60
\mathchardef\smallfrown="3\msx@61
\mathchardef\Subset="3\msx@62
\mathchardef\Supset="3\msx@63
\mathchardef\Cup="2\msx@64

\mathchardef\Cap="2\msx@65

\mathchardef\curlywedge="2\msx@66
\mathchardef\curlyvee="2\msx@67
\mathchardef\leftthreetimes="2\msx@68
\mathchardef\rightthreetimes="2\msx@69
\mathchardef\subseteqq="3\msx@6A
\mathchardef\supseteqq="3\msx@6B
\mathchardef\bumpeq="3\msx@6C
\mathchardef\Bumpeq="3\msx@6D
\mathchardef\lll="3\msx@6E

\mathchardef\ggg="3\msx@6F

\mathchardef\circledS="0\msx@73
\mathchardef\pitchfork="3\msx@74
\mathchardef\dotplus="2\msx@75
\mathchardef\backsim="3\msx@76
\mathchardef\backsimeq="3\msx@77
\mathchardef\complement="0\msx@7B
\mathchardef\intercal="2\msx@7C
\mathchardef\circledcirc="2\msx@7D
\mathchardef\circledast="2\msx@7E
\mathchardef\circleddash="2\msx@7F
\def\ulcorner{\delimiter"4\msx@70\msx@70 }
\def\urcorner{\delimiter"5\msx@71\msx@71 }
\def\llcorner{\delimiter"4\msx@78\msx@78 }
\def\lrcorner{\delimiter"5\msx@79\msx@79 }
\def\yen{\mathhexbox\msx@55 }
\def\checkmark{\mathhexbox\msx@58 }
\def\circledR{\mathhexbox\msx@72 }
\def\maltese{\mathhexbox\msx@7A }
\mathchardef\lvertneqq="3\msy@00
\mathchardef\gvertneqq="3\msy@01
\mathchardef\nleq="3\msy@02
\mathchardef\ngeq="3\msy@03
\mathchardef\nless="3\msy@04
\mathchardef\ngtr="3\msy@05
\mathchardef\nprec="3\msy@06
\mathchardef\nsucc="3\msy@07
\mathchardef\lneqq="3\msy@08
\mathchardef\gneqq="3\msy@09
\mathchardef\nleqslant="3\msy@0A
\mathchardef\ngeqslant="3\msy@0B
\mathchardef\lneq="3\msy@0C
\mathchardef\gneq="3\msy@0D
\mathchardef\npreceq="3\msy@0E
\mathchardef\nsucceq="3\msy@0F
\mathchardef\precnsim="3\msy@10
\mathchardef\succnsim="3\msy@11
\mathchardef\lnsim="3\msy@12
\mathchardef\gnsim="3\msy@13
\mathchardef\nleqq="3\msy@14
\mathchardef\ngeqq="3\msy@15
\mathchardef\precneqq="3\msy@16
\mathchardef\succneqq="3\msy@17
\mathchardef\precnapprox="3\msy@18
\mathchardef\succnapprox="3\msy@19
\mathchardef\lnapprox="3\msy@1A
\mathchardef\gnapprox="3\msy@1B
\mathchardef\nsim="3\msy@1C
\mathchardef\ncong="3\msy@1D

\mathchardef\varsubsetneq="3\msy@20
\mathchardef\varsupsetneq="3\msy@21
\mathchardef\nsubseteqq="3\msy@22
\mathchardef\nsupseteqq="3\msy@23
\mathchardef\subsetneqq="3\msy@24
\mathchardef\supsetneqq="3\msy@25
\mathchardef\varsubsetneqq="3\msy@26
\mathchardef\varsupsetneqq="3\msy@27
\mathchardef\subsetneq="3\msy@28
\mathchardef\supsetneq="3\msy@29
\mathchardef\nsubseteq="3\msy@2A
\mathchardef\nsupseteq="3\msy@2B
\mathchardef\nparallel="3\msy@2C
\mathchardef\nmid="3\msy@2D
\mathchardef\nshortmid="3\msy@2E
\mathchardef\nshortparallel="3\msy@2F
\mathchardef\nvdash="3\msy@30
\mathchardef\nVdash="3\msy@31
\mathchardef\nvDash="3\msy@32
\mathchardef\nVDash="3\msy@33
\mathchardef\ntrianglerighteq="3\msy@34
\mathchardef\ntrianglelefteq="3\msy@35
\mathchardef\ntriangleleft="3\msy@36
\mathchardef\ntriangleright="3\msy@37
\mathchardef\nleftarrow="3\msy@38
\mathchardef\nrightarrow="3\msy@39
\mathchardef\nLeftarrow="3\msy@3A
\mathchardef\nRightarrow="3\msy@3B
\mathchardef\nLeftrightarrow="3\msy@3C
\mathchardef\nleftrightarrow="3\msy@3D
\mathchardef\divideontimes="2\msy@3E
\mathchardef\varnothing="0\msy@3F
\mathchardef\nexists="0\msy@40
\mathchardef\mho="0\msy@66
\mathchardef\eth="0\msy@67
\mathchardef\eqsim="3\msy@68
\mathchardef\beth="0\msy@69
\mathchardef\gimel="0\msy@6A
\mathchardef\daleth="0\msy@6B
\mathchardef\lessdot="3\msy@6C
\mathchardef\gtrdot="3\msy@6D
\mathchardef\ltimes="2\msy@6E
\mathchardef\rtimes="2\msy@6F
\mathchardef\shortmid="3\msy@70
\mathchardef\shortparallel="3\msy@71
\mathchardef\smallsetminus="2\msy@72
\mathchardef\thicksim="3\msy@73
\mathchardef\thickapprox="3\msy@74
\mathchardef\approxeq="3\msy@75
\mathchardef\succapprox="3\msy@76
\mathchardef\precapprox="3\msy@77
\mathchardef\curvearrowleft="3\msy@78
\mathchardef\curvearrowright="3\msy@79
\mathchardef\digamma="0\msy@7A
\mathchardef\varkappa="0\msy@7B
\mathchardef\hslash="0\msy@7D
\mathchardef\hbar="0\msy@7E
\mathchardef\backepsilon="3\msy@7F
\def\Bbb{\ifmmode\let\next\Bbb@\else
 \def\next{\errmessage{Use \string\Bbb\space only in math mode}}\fi\next}
\def\Bbb@#1{{\Bbb@@{#1}}}
\def\Bbb@@#1{\fam\msyfam#1}

\catcode`\@=12

\newtheorem{thm}{Theorem}[section]
\def\bthm{\begin{thm}}
\def\ethm{\end{thm}}
\newtheorem{propn}[thm]{Proposition}
\def\bpropn{\begin{propn}}
\def\epropn{\end{propn}}
\newtheorem{cor}[thm]{Corollary}
\def\bcor{\begin{cor}}
\def\ecor{\end{cor}}
\newtheorem{lemma}[thm]{Lemma}
\def\blemma{\begin{lemma}}
\def\elemma{\end{lemma}}
\def\beqn{\begin{equation}}
\def\eeqn{\end{equation}}
\def\bproof{\noindent Proof - }
\def\eproof{\ $\Box$ \vskip .1in}

\newenvironment{example1}{\begin{example2} \rm}{\end{example2}}
\newtheorem{example2}[thm]{Example}

\newenvironment{remark1}{\begin{remark2} \rm}{\end{remark2}}
\newtheorem{remark2}[thm]{Remark}

\newenvironment{defn1}{\begin{defn2} \rm}{\end{defn2}}
\newtheorem{defn2}[thm]{Definition}
\def\bdefn{\begin{defn1}}
\def\edefn{\end{defn1}}
\def\beqn{\begin{equation}}
\def\eeqn{\end{equation}}
\def\bea*{\begin{eqnarray*}}
\def\eea*{\end{eqnarray*}}

\def\benum{\begin{enumerate}}
\def\eenum{\end{enumerate}}

\def\german{\frak}
\newcommand{\R}{{\Bbb R}}

\newcommand{\C}{{\Bbb C}}
\newcommand{\Z}{{\Bbb Z}}

\def\ra{\rightarrow}
\def\o{\sigma}
\def\a{\alpha}

\def\Del{\Delta}

\def\l{\langle}
\def\r{\rangle}
\def\om{\omega}
\def\lam{\lambda}
\def\pwc{{\rm int}\ \t_+^*}
\def\cpwc{\t_+^*}

\def\tPhi{\tilde \Phi}
\def\Phinv{\Phi^{-1}}

\def\g{\german g}
\def\t{\german t}

\def\u{\german u}
\def\ann{\rm ann}
\def\rank{{\rm rank}}

\def\exp{{\rm exp}}

\def\Tr{{\rm Tr}}

\def\id{{\rm id}}
\def\ol{\overline}

\def\Go{ G_\o }

\def\go{ \g_\o}

\def\Uo{ U_\o }

\def\zo{z_\o}
\def\Zo{Z_\o}
\def\ra{\rightarrow}

\begin{document}

\title{Multiplicity-free Hamiltonian actions need not be
K\"ahler}


\author{Chris Woodward
\thanks{Supported by an ONR Graduate Fellowship.
E-mail:  woodward@math.mit.edu.} \\
\small Dept. of Mathematics, Rm 2-229, M.I.T., Cambridge MA 02139}

\maketitle

\begin{abstract}
In this note we show that
Tolman's example (of a six dimensional Hamiltonian $T^2$-space
with isolated fixed points and no compatible K\"{a}hler structure)
can be constructed from the flag variety $U(3)/U(1)^3$ by
$U(2)$-equivariant symplectic surgery.
This implies that Tolman's space has a ``transversal multiplicity-free''
action of $U(2)$ and that Delzant's theorem
``every compact multiplicity-free torus action is K\"{a}hler''
\cite{D1} does not generalize to non-abelian actions.\footnote{I have learned
from F. Knop  that he has
independently constructed examples of non-K\"{a}hler multiplicity-free
actions, by a different method.}

Comments and criticisms of this preliminary version are welcome.
\end{abstract}

{\footnotesize
\tableofcontents}

\section{Introduction}

In February 1995 S. Tolman \cite{T} presented a non-K\"{a}hler Hamiltonian
action with isolated fixed points.  Here, non-K\"{a}hler
means that there is no invariant complex structure $J$ such that
the equation
$$ g(X,Y) = \om(X,JY), $$
where $\om$ is the symplectic form,
defines a Riemannian metric $g$.  The existence of such an example
answers the crucial question about Hamiltonian group actions
of whether a ``generic'' Hamiltonian action has a compatible
K\"{a}hler structure.
The use of the word ``generic'' here is vague.  Commonly,
``generic'' has been understood to mean that a maximal torus
should act with isolated fixed points.
The question above is crucial because the definition of the ``quantization''
of a Hamiltonian action, and the proofs of theorems such the
Guillemin-Sternberg multiplicity formula (see \cite{GSm}, \cite{M}),
are easier if one can choose a K\"ahler polarization.
Tolman's example implies that the class of examples studied
in equivariant symplectic geometry is much larger than
that studied in equivariant K\"ahler geometry, and that
K\"ahler techniques are not generally valid.

The main point of this note is to give a construction of Tolman's
example which shows that its symmetry group is larger than the
one originally proposed.  Namely, Tolman's example has a
``multiplicity-free'' action of $U(2)$.   Multiplicity-free
means that all of the symplectic reduced spaces associated to
the action are zero dimensional.
The most well-known examples of multiplicity-free
Hamiltonian actions are the Hamiltonian actions associated to toric
varieties and coadjoint orbits.  In fact, Delzant \cite {D1} proved that
any compact multiplicity-free torus action is the
Hamiltonian action associated to a smooth projective toric
variety.  The main result of this note implies that Delzant's
theorem does not generalize to non-abelian actions.  That is,
there can be only a limited relationship between
the theory of multiplicity-free actions in symplectic geometry
and the theory of spherical varieties (the non-abelian
generalization of toric varieties in algebraic geometry;
see \cite{LV},\cite{B},\cite{Kn}.)

I would like to thank S. Tolman for giving us a truly wonderful example,
and her, Y. Karshon, and V. Guillemin for suggesting improvements.

\section{Tolman's example}

Tolman \cite{T} proves that a symplectic gluing of two halves
of two six-dimensional Hamiltonian $T^2$-spaces, $M_1$ and $M_2$,
results in a non-K\"ahler Hamiltonian $T^2$-space $M_3$.
The gluing is represented by the following
picture, which will be explained in a moment.

\begin{center}
\setlength{\unitlength}{0.0075in}
\begin{picture}(680,195)(0,-10)
\put(20,140){\blacken\ellipse{6}{6}}
\put(20,60){\blacken\ellipse{6}{6}}
\put(60,20){\blacken\ellipse{6}{6}}
\put(140,20){\blacken\ellipse{6}{6}}
\put(140,60){\blacken\ellipse{6}{6}}
\put(60,140){\blacken\ellipse{6}{6}}
\put(20,140){\blacken\ellipse{6}{6}}
\path(20,140)(20,60)(60,20)
	(140,20)(140,60)(60,140)(20,140)
\path(60,140)(60,20)
\path(20,60)(140,60)
\path(20,140)(140,20)
\put(220,180){\blacken\ellipse{6}{6}}
\put(220,20){\blacken\ellipse{6}{6}}
\put(380,20){\blacken\ellipse{6}{6}}
\put(260,100){\blacken\ellipse{6}{6}}
\put(260,60){\blacken\ellipse{6}{6}}
\put(300,60){\blacken\ellipse{6}{6}}
\path(220,180)(220,20)(380,20)
	(220,180)(260,100)(260,60)
	(300,60)(260,100)
\path(260,60)(220,20)
\path(300,60)(380,20)
\path(420,100)(460,100)
\path(452.000,98.000)(460.000,100.000)(452.000,102.000)
\put(520,180){\blacken\ellipse{6}{6}}
\put(520,60){\blacken\ellipse{6}{6}}
\put(560,20){\blacken\ellipse{6}{6}}
\put(600,60){\blacken\ellipse{6}{6}}
\put(560,100){\blacken\ellipse{6}{6}}
\put(680,20){\blacken\ellipse{6}{6}}
\put(565,105){\makebox(0,0)[lb]{$p$}}
\put(565,25){\makebox(0,0)[lb]{$q$}}
\put(605,65){\makebox(0,0)[lb]{$r$}}
\path(520,180)(520,60)(560,20)
	(680,20)(520,180)
\path(520,180)(560,100)(560,20)
\path(520,60)(600,60)(680,20)
\path(600,60)(560,100)
\dashline{4.000}(0,140)(140,0)
\dashline{4.000}(200,140)(340,0)
\dashline{4.000}(500,140)(640,0)
\end{picture}

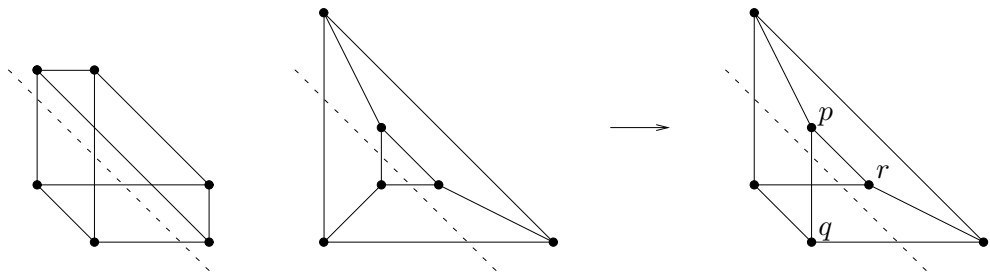
Figure 1:  Tolman's picture
\end{center}

The picture shows the images under the moment map of the connected
components of the non-principal orbit-type strata of $M_1$,
$M_2$, and $M_3$, respectively, from left to right.
That is, the points are the images of the
fixed points, and the line segments are the images of the submanifolds
fixed by circle subgroups of $T = T^2$.  (In general,
there are orbit-type strata with discrete isotropy groups, but
in these examples it turns out that there aren't any.)  Note that
an intersection of two lines is not necessarily the image of a fixed
point.

$M_1$ is a generic coadjoint orbit of $U(3)$, with the Hamiltonian
action of $U(3)$ restricted to $T = U(1)^2 \times {\rm id} \subset
U(3)$.  $M_2$ is a toric variety, with the action
restricted from $U(1)^3$ to $T$.  The polytope $P$ associated
to this toric variety is obtained by making the vertices of the
outer triangle have $z$-coordinate $0$ and those of the inner triangle
$z$-coordinate $1$, and giving each vertex $x$ and $y$ coordinates as
drawn.  The picture above is then the projection of $P$ onto the
$x$-$y$ plane.\footnote{This toric variety is the blow-up of
a weighted projective space at a fixed point.}
The pictures are representations of an invariant which Tolman
calls the X--ray.

To describe this invariant, let $(M,\om)$
be a Hamiltonian $T$-space with moment map $\Phi:M \ra \t^*$,
and let $\chi$ be the set of connected components of orbit-type strata.
(Recall that the orbit-type stratum $M_H$ corresponding to a
subgroup $H \subset T$ is the set of points $m \in M$ such that
the isotropy subgroup $T_m$ equals $H$.)
\bdefn  The X--ray of M is the set of images
$\{ \Phi(X) \ \vert \ X \in \chi \}$
partially ordered by the relation $\{ ( \Phi(X_1), \Phi(X_2) )
\ \vert \ \ X_1 \subset \ol{X_2},\ \  X_1,X_2 \in \chi \}$.
\edefn

\subsection{X--ray computations}

Although the geometry of the spaces $M_1$ and $M_2$ is well-known,
for the convenience of the reader we will sketch the computations
of their X--rays.   To compute the X--ray for $M_1$,
we can use the general fact that

\blemma \label{LineLemma}
If $M$ is a Hamiltonian $T$-space with moment map $\Phi:M \ra \t^*$,
and if $\a \in \t^*$ is a weight with multiplicity one
for the action of $T$ on the tangent space of $M$
at a fixed point $m$, then
there is an $X \in \chi$ such that $\Phi(X)$ lies in
$\Phi(m) + \R_+ \a$.
\elemma

\bproof The equivariant Darboux theorem produces
a local isomorphism $\varphi: (T_m M, \om_m) \ra (M, \om) $.
$X$ is locally just the image of the weight space $(T_m M)_\a$ under
$\varphi$.
The moment map for the action of $T$ on $(T_m M)_\a
\cong \C$ is $ z \mapsto \Phi(m) + \a \vert z \vert^2/2  $,
which proves the lemma.
\eproof
Suppose that $M_1 = U(3) \lam$ for some generic $\lam \in \u(3)^*$.
(We could also take $M_1$ to be an $SU(3)$-coadjoint orbit, and restrict
to the action of the maximal torus.)
The $T$-fixed points are the $U(1)^3$-fixed points, that is,
the elements of $W \lam$, where
$W$ is the Weyl group of $U(1)^3 \subset U(3)$.
The weights of $T$ at a fixed point
$w \lam $ are $\pm \a_1, \pm \a_2 $ and $\pm \a_3$, where $\a_1$,
$\a_2$, and $\a_3$ are the positive weights restricted to
$\t \subset \u(1)^3$.   With respect to standard bases,
the positive roots are
$$ (1,-1,0), \ \ (0,1,-1), \  \ (1,0,-1) $$
and so
$$ \a_1 = (1,-1) , \ \ \a_2 = (0,1), \ \ \a_3 = (1,0). $$
Therefore in a neighborhood of each
$w \lam $ there are three non-principal strata
fixed by circle subgroups of $T$.  Since any two of $\a_1,\a_2$,
and $\a_3$ are a basis of the weight lattice, there are
no strata with discrete isotropy groups in these neighborhoods.
If $X \in \chi$ has isotropy group $H$, then
$\ol{X}$ is a component of the fixed point set of $H$, that is,
a compact submanifold, and therefore $\ol{X}$ must contain $T$-fixed points.
Thus, the strata which meet the local models around the
fixed points are the only ones.
Furthermore, this implies that the endpoints of any line segment
$\Phi(\ol{X}),\ X \in \chi$ must be images of fixed points.
It follows that the X--ray of $M_1$ is as shown in Figure 1.

To compute that X--ray for $M_2$, let $\Phi_3:M_1 \ra (\t^3)^*$
be the moment map for the 3-torus, $P = \Phi_3(M)$ the moment polytope,
and $F$ a face of $P$.  The following is a general fact about
toric varieties.  For a symplectic proof, see \cite{D1}.
\blemma \label{DelLem}
The isotropy subgroup of any point $m$ in $\Phinv_3(F)$ is
the subgroup $T_F^3$ with Lie algebra $\ann(F)$,
where $\ann(F) \subset \t^3$ is the annihilator of $F$.
\elemma
The isotropy subgroup for the action of the 2-torus $T$ on
$\Phinv_3(F)$ is therefore $T \cap T_F^3$.  Furthermore, the
moment map $\Phi_2$ for $T$ is $\Phi_3$ composed with projection onto $\t^*$.
If $F$ is a $0$ or $1$-dimensional face of $P$, it follows that
$\Phi_3^{-1}(F)$
is a component of a non-principal orbit-type stratum of the $T^2$ action
whose whose image under $\Phi_2$ is the projection of
$F$ onto $\t^*$.  If $F$ is a
$2$-dimensional face, $T \cap T_F^3 = \{ {\rm id } \} $, so $\Phinv_3(F)$
is part of the principal orbit-type stratum for the $T$-action.
This implies that the X--ray is as shown in Figure 1.

\subsection{Non-existence of a compatible K\"ahler structure}

In \cite{T} Tolman (using a simple case of the classification of
Hamiltonian actions on $4$-manifolds by Audin-Ahara-Hattori-Karshon)
shows that neighborhoods in $M_1$ and $M_2$ of the inverse images of
the dotted lines shown in Figure 1 are equivariantly symplectomorphic.
Tolman's space $M_3$ can then be constructed by gluing one half of $M_1$
together with one half of $M_2$.

Suppose that $M_3$ is K\"ahler.  In this situation,
the action complexifies and a theorem of Atiyah \cite{A}
describes the images of the orbits of the complex torus:
\bthm[Atiyah]  Let $M$ be a compact K\"{a}hler manifold with a Hamiltonian
action of a torus $T$.  Let $T_\C$ be the complexified torus
and $Y$ an orbit of $T_\C$.  Then
\begin{itemize}
\item[(1)]  $\Phi(\ol{Y})$ is a convex polytope $\Del$.
\item[(2)]  If $F \subset \Del$ is any face, then
$\Phinv(F) \cap \ol{Y}$ is a single $T_\C$ orbit, $Y_F$.
\item[(3)]  $\Phi$ induces a homeomorphism $\ol{Y}/T \cong \Del$.
\end{itemize}
\ethm
Let $C$ be the cone based at $p$ in Figure 1, and generated by the line
segments $\ol{pq}$ and $\ol{pr}$.  Tolman notes that there is a
orbit $Y$ of $T_\C$ with the property that $\Phi(\ol{Y}) = C$ near $p$.
$Y$ can be constructed using an equivariant local chart for the action
of $T_\C$.
An alternative construction goes as follows.  Let $U$ be the unstable
manifold of $m$ with respect to a Morse function of the form
$\l \Phi, v \r$, for a generic $v \in \t$ such that $(v,q-p) < 0$
and $(v,r-p) < 0$, so that $U$ is $4$-dimensional.
Since the gradient flow is given by the
action of a one-parameter subgroup of $T_\C$,
$U$ is a $T_\C$-space and splits into $T_\C$-orbits.
Since these orbits are even real-dimensional
there is one orbit which is dense in $U$, which we define to be $Y$.
If $X_{pq}$ and $X_{pr} \in \chi$ are such that
$\Phi(\ol{X_{pq}}) = \ol{pq}$ and $\Phi(\ol{X_{pr}}) = \ol{pr}$ then $U$
must contain $X_{pq}$ and $X_{qr}$, since the limit point for the
gradient flow on these submanifolds must be a $T$-fixed point and
the only possibility is $m$.
By the application of Atiyah's theorem, $\Phi(\ol{Y})$
is a convex polytope $\Del$.  The submanifolds $Y_F$ are
Hamiltonian $T$-spaces, and since $\Phi(Y_F) = F$ we have
by definition of the moment map that $Y_F$ is fixed by the
subgroup of $T$ with Lie algebra $\ann(F)$.  In particular,
the set $\{ p \in \ol{Y} \vert \dim T_p > 0 \}$ of points
with positive dimensional isotropy subgroups must map to the boundary
$\partial \Del $ of $\Del$.
Therefore $\ol{pq} \cup \ol{pr} \subset \partial \Del$,
so that $\Del = C$ locally, as claimed.

The only images of fixed points in the cone $C$ are $p,q$, and $r$
and so $\Del$ must be the convex hull of $p,q$ and $r$.  The line
segment $\ol{qr}$ minus its endpoints is therefore a face
$F$ of $\Del$.  $Y_F$ is fixed by a circle subgroup,
and therefore must be contained in some $X \in \chi$ such that
$\Phi(\ol{X})$ contains $F$.  But from Figure 1 we see
that there is no such $X$, which is a contradiction.\footnote{One
can give an alternative proof that $M_3$ is not K\"ahler using
Witten's equivariant holomorphic Morse inequalities, in a form
due to Siye Wu.  This observation arose during a discussion with
Wu and Tolman, whom we thank.}

This argument is a particular example of a more general criterion
which Tolman shows is necessary for the existence of a K\"ahler
structure, in terms of the X--ray.  We will not discuss the more general
statement here, since we want to focus on this example.

\section{Construction by $U(2)$-equivariant symplectic surgery}

In this section we give an alternative construction of Tolman's
example which uses the ``symplectic cutting'' technique of E. Lerman \cite{L}.
The advantage is that
(1) the construction is more explicit, and (2) the construction
shows that the example has a ``transversal'' multiplicity-free $U(2)$-action,
which explains the $\Z/2$ symmetry in Tolman's picture.   The symmetry
comes from the action of the Weyl group $\Z/2$ of $U(2)$.  Transversal
means that the moment map is transversal to the Cartan subalgebra.  This
implies that the action fits into the classification of \cite{W}.

\subsection{Lerman's definition of symplectic cutting.}

Let $(N,\om_N)$ be a Hamiltonian $G$-space and $\mu:N \ra \R$ the
moment map for a $G$-equivariant $S^1$ action on $N$.  Let $a \in \R$
be a regular value of $\mu$, let $N_a = \mu^{-1}(a)$ be the
reduced space and let $N_{< a}$ be the subset $\mu^{-1}(-\infty,a)
\subset N$.

\blemma[Lerman]
The union $N_{\le a} = N_a \cup N_{< a}$ has the structure
of a Hamiltonian $G$-orbifold.
Furthermore, if $N_a$ is smooth then $N_{\le a} $ is smooth also.
\elemma

\bproof   Let $N \times \C$
be the product with symplectic structure $ \pi_1^* \om_N + \pi_2^*
(dz \wedge  d\ol{z}) /2i  $, where $\pi_1$ and $\pi_2$ are the projections.
Define $\nu : N \times \C \ra \R$ by
$$ \nu(n,z) = \mu(n) + \vert z \vert^2 /2 $$
so that $\nu$ is the moment map of the diagonal action of $S^1$ on
$N \times \C$,
which is equivariant with respect to the action of $G$ on the left factor.
Let $N_{\le a}$ be the reduction of $N \times \C$ at $a$.  Then
we can write
$$N_{\le a} \cong \mu^{-1}(a)/S^1 \cup \mu^{-1}(-\infty,a)$$
as claimed.
\eproof
The space $N_{\le a}$ is called the symplectic cut of $N$ at $a$.
It is easy to check that the identification of a dense subset of $N_{\le a}$
with $N_{< a} \subset N$ is an equivariant symplectomorphism.
This implies that
$N_{\le a}$ is defined even if $\mu$ is only smooth in a neighborhood
$U$ of $\mu^{-1}(a)$.  (That is, $\mu$ only defines an $S^1$ action locally.)
Indeed, we can assume that $U = \mu^{-1}(b,c)$.
The symplectic cut $U_{\le a}$ is well-defined, and
we define (see \cite{L})
\bdefn  Let $N_{\le a}$ be the union of $U_{\le a}$ and
$N_{< a}$ glued together over $U_{< a}$.
\edefn

\subsection{Equivariant symplectic surgery}

If $M$ is a Hamiltonian $G$-space, there is a natural set of
(locally defined) $G$-equivariant circle actions which we can
use for symplectic cutting.  This follows from a version of the
Guillemin-Sternberg symplectic cross-section theorem
\cite[Theorem 26.2]{GSstp}.  Let $\cpwc$ be a closed positive
Weyl chamber.  For each Weyl wall $\o \subset \cpwc$
(not necessarily codimension 1) let $\Go$ be the stabilizer
of any point in $\o$.
\bthm[Guillemin-Sternberg]
Let $\o \subset \cpwc$ be a Weyl wall,
and let $\Uo \subset \go^*$ be the maximal set such that $x \in \Uo$
implies $G_x \subset \Go$.  Then
\benum
\item  $\Phinv(\Uo)$ is a Hamiltonian $\Go$-space, called the
symplectic cross-section for $\o$.   Note that
$\Phinv(G \Uo) \cong G \times_{\Go} \Phinv(\Uo) $.
\item  If $\Zo$ is the center of $\Go$, then we can define
a {\em new} action of $\Zo$ by extending the action
on $ \Phinv(G \Uo)$ by $G$-equivariance.
In particular if $\o = \pwc$ then we have a new $G$-equivariant
action of $T$ on the (dense if non-empty) subset
$G \Phinv(\Uo) = \Phinv(\g^*_{\rm reg})$.
\eenum
\ethm
We call the new action of $\Zo$ the right action, and denote
it by $\rho$.  The justification for this terminology is that if
$m \in \Uo$ and $g \in G$ then the new action is given by
$$ \rho(z) ( gm) = g (\rho(z) m) = g (zm) = (gz)m $$
and therefore the isomorphism  of the orbit $Gm$ with $G/G_m$
identifies $\rho$ with the right action of $\Zo$ on $G/G_m$
(which is well-defined since $G_m$ is contained in $\Go$
and so commutes with $\Zo$.)

Let $q: \g^* \ra \cpwc$ be the quotient map, and define
$\tPhi = q \circ \Phi$.
The following proposition goes back in some form to Thimm
(see e.g. \cite{GSgc}.)
\bpropn  The composition of
$\tPhi$ with the projection $\pi_\o:\t^* \ra \zo^*$ is a moment
map for the right action of $\Zo$.
\epropn
For convenience, we recall the proof.
Let $X \in \zo$, and let $X_L^\#$ and $X_R^\#$ be the generating vector
fields of the left and right actions of the one-parameter subgroup
$S^1 = \exp(t X)$.
Since $G \Uo = G \times_{\Go} \Uo$, and
$\l \tPhi, X \r = \l \Phi, X \r$ on $\Uo$, we have that
$\l \tPhi, X \r$ is smooth on $G \Uo$.
Let $Y$ the Hamiltonian vector field of $\l \tPhi, X \r$.
We must show $X_R^\# = Y$.  Since both are $G$-invariant
it suffices to show this at at $\Phinv(\Uo)$, where
$X_R^\# = X_L^\#$.  Since $\tPhi$ is $G$-invariant,
$Y \Phi = 0$, and so $Y$ is tangent to $\Phinv(\Uo)$.
If $j$ is the inclusion of $\Phinv(\Uo)$ in $M$ then
$$ \imath(Y) (j^* \om)  = d  \l j^* \tPhi, X \r = d \l j^* \Phi, X \r  =
\imath(X_L^\#) (j^* \om). $$
Since $\Phinv(\Uo)$ is a symplectic submanifold, this implies
that $Y = X^\#_L$ on $\Phinv(\Uo)$
as required.
\eproof

Since symplectic cutting is local, if $a \in \R$ is such that
$\mu = \l \tPhi, X \r$ is smooth at $\mu^{-1}(a)$, and if
furthermore the right action of $S^1 = \exp(t X) \subset T$
is free on $\mu^{-1}(a)$,
then the symplectic cut of $M$ at $a$ is a Hamiltonian $G$-space
$M_{\le a}$ with moment
polytope $$\Del_{\le a} =  \{ v \in \Del \ \vert \ \l v , X \r \le a \}. $$

Note that
by equivariance, in order to check that the right action of $S^1$ is free
is suffices to check that the left action of $S^1$ is free
on $\Phinv(\cpwc) \cap \mu^{-1}(a)$.   Also, the condition that
$\mu$ be smooth at $\mu^{-1}(a)$ is equivalent to requiring that
the hyperplane $H = \{ v \in \t^* \ \vert \ \l v , X \r = a \}$
meets perpendicularly every Weyl wall $\o$ such that
$\o \cap H \cap \Del$ is non-empty.

\subsection{Construction of Tolman's example}

Let $M_1 = U(3)/U(1)^3$ and consider the Hamiltonian action of
$U(2) \subset U(3)$ obtained by restriction.  Here $U(2)$ is
embedded in $U(3)$ by the map $ A \ra {\rm diag}(A,1)$.
Let $\Phi':M_1 \ra \u(2)^*$
be the moment map, $T \subset U(2)$ the diagonal maximal torus,
and $\Del' \subset \cpwc$ the moment polytope of $M_1$
pictured on the upper left.
The other three pictures will be explained later.

\begin{center}
\setlength{\unitlength}{0.0075in}
\begin{picture}(500,475)(0,-10)
\path(340,0)(340,40)(380,40)
	(460,0)(340,0)(300,40)
	(300,160)(340,80)(340,40)(300,40)
\put(340,0){\blacken\ellipse{6}{6}}
\put(340,80){\blacken\ellipse{6}{6}}
\put(380,40){\blacken\ellipse{6}{6}}
\put(460,0){\blacken\ellipse{6}{6}}
\put(300,40){\blacken\ellipse{6}{6}}
\put(300,160){\blacken\ellipse{6}{6}}
\put(340,5){\makebox(0,0)[lb]{$x_1$}}
\put(380,45){\makebox(0,0)[lb]{$x_2$}}
\put(460,5){\makebox(0,0)[lb]{$x_3$}}
\put(300,45){\makebox(0,0)[lb]{$wx_1$}}
\put(340,85){\makebox(0,0)[lb]{$wx_2$}}
\put(300,160,5){\makebox(0,0)[lb]{$wx_3$}}
\path(340,80)(380,40)
\path(300,160)(460,0)
\path(4.243,267.071)(0.000,260.000)(7.071,264.243)
\path(0,260)(200,460)
\path(195.757,452.929)(200.000,460.000)(192.929,455.757)
\path(200,260)(200,300)(40,300)
	(40,260)(200,260)
\put(200,260){\blacken\ellipse{6}{6}}
\put(200,300){\blacken\ellipse{6}{6}}
\put(40,260){\blacken\ellipse{6}{6}}
\dashline{4.000}(200,240)(0,340)
\path(304.243,267.071)(300.000,260.000)(307.071,264.243)
\path(300,260)(500,460)
\path(495.757,452.929)(500.000,460.000)(492.929,455.757)
\path(340,300)(340,260)(460,260)
	(380,300)(340,300)
\put(340,260){\blacken\ellipse{6}{6}}
\put(380,300){\blacken\ellipse{6}{6}}
\put(460,260){\blacken\ellipse{6}{6}}
\put(340,265){\makebox(0,0)[lb]{$x_1$}}
\put(380,305){\makebox(0,0)[lb]{$x_2$}}
\put(460,265){\makebox(0,0)[lb]{$x_3$}}
\path(240,360)(300,360)
\path(292.000,358.000)(300.000,360.000)(292.000,362.000)
\path(40,200)(0,200)(0,40)(40,40)
\path(40,40)(40,200)
\path(40,40)(40,0)(200,0)
	(200,40)(40,40)
\path(0,40)(40,0)
\path(40,200)(200,40)
\path(0,200)(200,0)
\path(220,100)(280,100)
\put(40,200){\blacken\ellipse{6}{6}}
\put(0,200){\blacken\ellipse{6}{6}}
\put(0,40){\blacken\ellipse{6}{6}}
\put(40,0){\blacken\ellipse{6}{6}}
\put(200,0){\blacken\ellipse{6}{6}}
\put(200,40){\blacken\ellipse{6}{6}}
\path(272.000,98.000)(280.000,100.000)(272.000,102.000)
\end{picture}

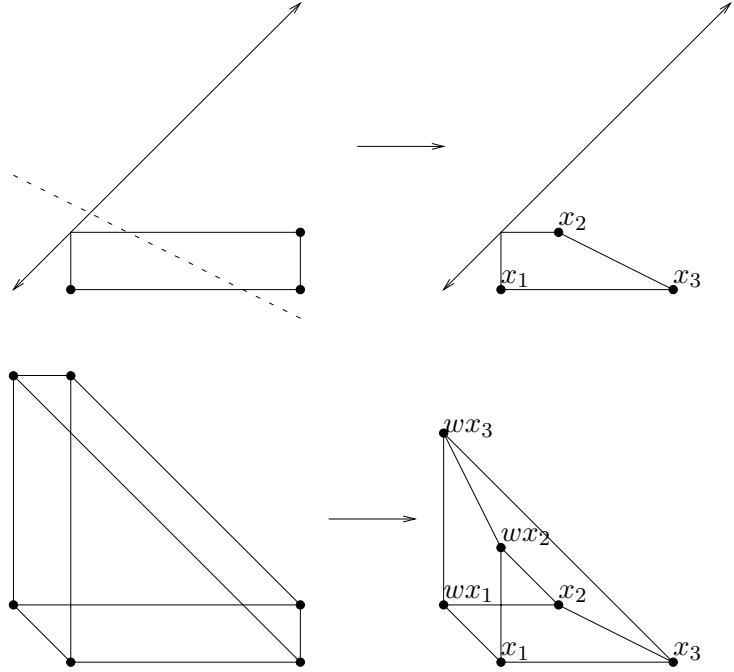
Figure 3:  Construction by symplectic cutting
\end{center}

$\Del'$ can be computed explicitly
(see \cite[page 364]{GSstp}) or as follows.\footnote{Note that this method
for computing the moment polytope only works in certain cases.}
The vertices of $\Del'$ contained in $\pwc$ are
images of $T$-fixed points, and there are three
elements of $W \lam$ which lie in the positive Weyl chamber
of the Weyl group of $U(2)$.
To find the vertices of $\Del'$
contained in the Weyl wall $\o$, note that if $F$ is any $1$-face of $\Del'$
not lying entirely in the $\o$
then $Y_F = \Phinv(\ol{F}) \cap Y_+$ is submanifold
fixed by a circle subgroup of $T$.\footnote{It follows from Kirwan's theorem
that
$Y_F$ is connected.}  One endpoint
of $F$ must lie in the positive Weyl chamber, so that
$Y_F$ must contain a fixed point $m$.
The tangent space $T_m M_1$ splits equivariantly
\beqn \label{split}
T_mM_1 \cong T_m Y_+ \oplus \u(2)^*/\t^*
\eeqn
and the weight of $T$ on $\u(2)^*/\t^*$
is $- \a_1$ where $\a_1$ is the positive root of $\u(2)$.
Suppose without loss of generality that $\Phi(m)$ is the lower
left vertex of $\Del'$, so that the weights of $T$ on $T_m M_1$
are $\{  -\a_1,\a_2,\a_3 \}$.  Equation (\ref{split}) implies
that weights of $T$ on $T_m Y_+$ are $\a_2$ and $\a_3$.
Since $T_m Y_F \subset T_m Y_+ $ is a weight space,
$F$ must be contained
in $\Phi(m) + \R_+ \a_2$ or $\Phi(m) + \R_+ \a_3$, but the
latter ray does not intersect $\o$.  This implies that
the additional vertex of $\Del'$ is the intersection of
$\Phi(m) + \R_+ \a_2$ with $\o$.

The map $\tPhi = q \circ \Phi'$ has two components, $\tPhi_1$ and
$\tPhi_2$.  Let $ \mu $ be the function
$\tPhi_1 + 2 \tPhi_2 $ and $S^1$ the circle subgroup
$\{ (z , z^2) \vert z \in U(1) \} \subset U(1)^2$
whose right action has moment map $\mu$.
Let $a \in \R$ be such that $\mu^{-1}(a)$
is the inverse image under $\tPhi$ of the dotted line shown in Figure 3.
Then $\mu$ is smooth at $\mu^{-1}(a)$, since $\mu^{-1}(a)$ lies entirely in
$G Y_+$.  The symplectic cut $(M_1)_{\le a}$ of $M_1$ at $a$ using
$\mu$ is therefore well-defined.

\bpropn  $(M_1)_{\le a}$ is a smooth Hamiltonian $U(2)$-space, whose
X--ray as a $T^2$-space is the same as in Tolman's
example.\footnote{Conjecturally, $(M_1)_{\le a}$ is equivariantly
$T^2$-symplectomorphic
to Tolman's space $M_3$.  This would follow from a theory of Hamiltonian
2-torus actions on 6-manifolds, which is work in progress by
S. Tolman and Y. Karshon.}
\epropn

To see that $(M_1)_{\le a}$ is smooth,
it's enough to check that the right action of $S^1$ is free on $\mu^{-1}(a)$,
or equivalently, since the right action is  $U(2)$-equivariant,
that $ S^1$ acts freely on $\mu^{-1}(a) \cap \Phinv(\cpwc)$.
This is shown either explicitly
(see \cite{W} page 6 or \cite{GSstp} page 364) or
using (slight extension of) Delzant's lemma \ref{DelLem}.
(See e.g. \cite[Lemma 4.2]{W}.)
\blemma[Delzant]
Let $M$ a multiplicity-free compact Hamiltonian $G$-space
with moment polytope $\Del$, let $F$ be a face of $\Del$
contained in $\pwc$, and let $m \in \Phinv(F)$.   Then
the isotropy subgroup $G_m$
of $m$ is connected and its Lie algebra is $\ann(F) \subset \t$.
\elemma
By the lemma the isotropy subgroups of points in
$\mu^{-1}(a) \cap Y_+ $ are $\id \times U(1)$ and the
trivial group.  These intersect $S^1$ trivially,
so $S^1$ acts freely on $\mu^{-1}(a)$.
Therefore, $(M_1)_{\le a}$ is a smooth Hamiltonian $U(2)$-space which
has the moment polytope shown on the above right.
Since $T$ acts freely on a dense subset of $Y_+$,
$U(2)$ acts freely on a dense subset of $GY_+$.  Since
$\dim((M_1)_{\le a}) = 6 = (\dim + \rank)U(2) $ we have
that $(M_1)_{\le a}$ is multiplicity-free.

Now consider the actions of $T = U(1)^2 \subset U(2)$ on $(M_1)_{\le a}$,
and let $\Phi:(M_1)_{\le a} \ra \t^*$ be the moment map.
Let us compute the X--ray for $(M_1)_{\le a}$.
Since $(M_1)_{\le a}$ is a symplectic cut, the $T$-fixed points are those
lying in $\mu^{-1}(-\infty,a)$ and the ``new'' fixed points
in $\mu^{-1}(a)/S^1$.  From Figure 3, we see that there is only
one ``old'' fixed point $q_1$ with $\Phi(q_1) = x_1$, and
by Delzant's lemma (or explicitly)
there are $T$-fixed points $q_2,q_3$ whose images under $\Phi$ are
$x_2$ and $x_3$ resp.  The points $wq_i, i = 1,2,3$, where
$w $ is the non-trivial element of the Weyl group, are also $T$-fixed points,
whose images under $\Phi$ are $wx_i$.

The lines in the X--ray are using lemma \ref{LineLemma}.
The splitting of $T_{q_i} (M_1)_{\le a}$ as in equation (\ref{split})
implies that
there are two weights of $T$ on $T_{q_i} (M_1)_{\le a}$ which are weights
of $T$ acting on $T_{q_i} (Y_+)_{\le a}$
and a third weight equal to  $ - \a_1$.
By lemma \ref{LineLemma} and Delzant's lemma
the two weights of $T_{q_i} (Y_+)_{\le a}$
are proportional to the  directions of the two edges of the moment
polytope at $\tPhi(q_i)$.
The weights at the $T$-fixed points $wq_i$ are obtained by Weyl reflection.
By the reasoning similar to the computation of the X--ray of $M_1$,
the X--ray can be computed from this local data,
and coincides with the X--ray of Tolman's example.
\eproof

\section{How many Hamiltonian actions are K\"ahler?}

One way of giving
meaning to this question is to consider a family of Hamiltonian
actions described by some parameter, and to describe for which
values of the parameter an invariant, compatible K\"ahler structure exists.
In this section, we describe a family of examples parametrized by the
integers.  Call a multiplicity-free
Hamiltonian $U(2)$-space a generalized Hirzebruch space
if its moment polytope is of the following form for some $n \in \Z$:

\begin{center}

\setlength{\unitlength}{0.0075in}
\begin{picture}(200,215)(0,-10)
\path(195.757,192.929)(200.000,200.000)(192.929,195.757)
\path(200,200)(0,0)
\path(4.243,7.071)(0.000,0.000)(7.071,4.243)
\texture{c0c0c0c0 0 0 0 0 0 0 0
	c0c0c0c0 0 0 0 0 0 0 0
	c0c0c0c0 0 0 0 0 0 0 0
	c0c0c0c0 0 0 0 0 0 0 0 }
\shade\path(60,60)(60,20)(200,20)
	(80,60)(60,60)
\path(60,60)(60,20)(200,20)
	(80,60)(60,60)
\put(150,40){\makebox(0,0)[lb]{slope $-1/n$}}
\end{picture}

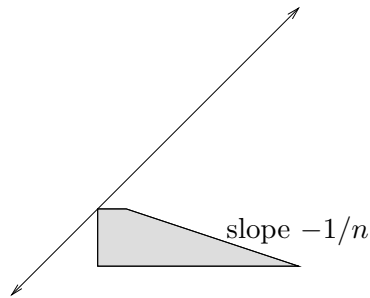
Figure 4:  Hirzebruch polytopes

\end{center}

The reason for this language is that these are the polytopes
of the famous Hirzebruch surfaces considered as toric varieties.
Generalized Hirzebruch spaces can be constructed
by symplectic cutting using the function
$\tPhi_1 + n \tPhi_2 $ on the flag variety $U(3)/U(1)^3$
as described above.  It follows from a result of Delzant \cite{D2}
(see also \cite{W}) that there is only one such space whose moment polytope
has slope $1/n$.  We denote this space by $H_n$.

For $n \ge 2$, Tolman's criterion shows $H_n$
cannot have a compatible, invariant K\"{a}hler structure.
The question of whether $H_n$ has a non-invariant compatible
K\"{a}hler structures is open.  There
is no known obstruction to the existence of a non-invariant
K\"{a}hler structure.  In particular, their cohomology rings
(which can be computed by studying the variation of the moment polytope
with respect to the cohomology class of the symplectic form)
have the hard Lefshetz property.

For $n=1$, the
function $\tPhi_1 + \tPhi_2 = \Tr(\Phi) = \Phi_{11} + \Phi_{22}$
is the Hamiltonian of the action of the center $Z \subset U(2)$,
which is an everywhere-defined action which preserves the K\"{a}hler
structure.  Therefore in this case the symplectic cuts {\em are}
K\"{a}hler, by \cite[Theorem 3.5]{GSm}.
This example is discussed in greater detail in \cite{W}.

For $n=0$ we get the flag variety $U(3)/U(1)^3$, considered
as Hamiltonian $U(2)$-space.

In the case $n<0$ F. Knop \cite{Knpc} has observed that there exist
compatible, invariant K\"ahler structures.  Knop uses
the theory of spherical varieties (see \cite{LV},\cite{Kn})
to construct equivariant
completions of the homogeneous space $GL(2,\C)/\C^*$.  The
existence of projective structures on these varieties
follows from work of Brion \cite[page 413]{Br}.
The proof will appear elsewhere.

To summarize, the results of this note together with Knop's
observation imply that

\bthm  If $n \ge 2$, then $H_n$ has no $T$-invariant
compatible K\"ahler structure.  Otherwise, $H_n$ is a K\"ahler
Hamiltonian $U(2)$-space.
\ethm

\end{document}